\documentclass[a4paper,11pt,oneside]{article}

\usepackage{amsmath,amssymb}
\usepackage{epsfig}
\usepackage[utf8]{inputenc}
\usepackage{textcomp}
\usepackage{algorithm}
\usepackage{algpseudocode}
\usepackage{color}
\usepackage{graphicx}
\usepackage{caption} 
\usepackage{subcaption}
\usepackage{bm}
\usepackage{amsthm}
\usepackage{mathrsfs}
\usepackage{dsfont}
\usepackage{url}

\usepackage{geometry}
\newtheorem{thm}{Theorem}[section]
\newtheorem{cor}[thm]{Corollary}
\newtheorem{definition}[thm]{Definition}

\newtheorem{proper}[thm]{Property}

\newtheorem{remark}[thm]{Remark}

\newcommand{\Prob}{\mathbb{P}}

\newcommand{\E}{\mathbb{E}}

\makeatletter
\def\blfootnote{\xdef\@thefnmark{}\@footnotetext}
\makeatother

\title{Optimal piecewise linear data compression for solutions of parametrized partial differential equations}

\author{Thomas DANIEL${}^{*,\dagger}$, Fabien CASENAVE${}^{*}$, Nissrine AKKARI${}^{*}$, \\ David RYCKELYNCK${}^{\dagger}$}

\begin{document}

\maketitle
\blfootnote{${}^{*}$ SafranTech, Rue des Jeunes Bois, Ch\^ateaufort, 78114 Magny-les-Hameaux (France).}
\blfootnote{${}^{\dagger}$ MINES ParisTech, PSL University, Centre des mat\'eriaux (CMAT), CNRS UMR 7633, BP 87, 91003 Evry (France).}
\blfootnote{Corresponding authors: \url{thomas.daniel@mines-paristech.fr} and \url{fabien.casenave@safrangroup.com}.}

\begin{abstract}
Model order reduction has been extensively studied over the last two decades. Projection-based methods such as the Proper Orthogonal Decomposition and the Reduced Basis Method enjoy the important advantages of Galerkin methods in the derivation of the reduced problem, but are limited to linear data compression for which the reduced solution is sought as a linear combination of spatial modes. Nonlinear data compression must be used when the solution manifold is not embedded in a low-dimensional subspace. Early methods involve piecewise linear data compression, by constructing a dictionary of reduced-order models tailored to a partition of the solution manifold. In this work, we introduce the concept of optimal partition of the solution manifold in terms of normalized Kolmogorov widths, and prove that the optimal partitions can be found by means of a representative-based clustering algorithm using the sine dissimilarity measure on the solution manifold.
\end{abstract}

\section{Introduction}
\label{sec:intro}

Let $\Omega$ be an open set in $\mathbb{R}^d$, and let us introduce a generic parametrized partial differential equation (pPDE):
\begin{equation}
\mathcal{D}(u;x) = 0,
\label{PDE}
\end{equation}
where $u\in \mathcal{H}\subset L^2 (\Omega)$ and $x\in\mathcal{X}\subset\mathbb{R}^p$ denote respectively the solution and the parameter, with $\mathcal{H}$ being a Hilbert space and $\mathcal{X}$ the parameter domain. For time-dependent problems, the time is considered as a parameter. We denote by $\|\cdot\|_{\mathcal{H}}$ the norm induced by the inner product $\langle \cdot , \cdot \rangle_{\mathcal{H}}$ of $\mathcal{H}$. The pPDE~\eqref{PDE} involves differential operators and initial conditions and/or boundary conditions, and is assumed to be well-posed in the sense of Hadamard: $\mathcal{X}\ni x \mapsto u \in \mathcal{H}$ solution of \eqref{PDE} is a continuous application, called the solution application. For clarity of presentation, this application is still denoted by $u$: for all $x\in\mathcal{X}$, $u(x)$ is the unique solution of~\eqref{PDE}. The solution manifold $\mathcal{M}$ is defined as the image of the solution application: $\mathcal{M} = u(\mathcal{X})$.
We suppose that $\mathcal{M}$ does not contain the zero solution: $\exists \ \lambda>0$ such that $\displaystyle \underset{x\in\mathcal{X}}{\rm inf}\|u(x)\|_{\mathcal{H}}=\underset{u\in\mathcal{M}}{\rm inf}\|u\|_{\mathcal{H}}>\lambda$.

Uncertainty quantification or optimization tasks involving the pPDE~\eqref{PDE} require computing approximate solutions for many values of the parameter $x$. These many-query tasks can become intractable if a high-dimensional approximation of~\eqref{PDE}, called High-Dimensional Model (HDM), is considered. The dimension of the HDM is denoted by $\mathcal{N}$. In such cases, reduced-order modeling can be used to derive a low-dimensional approximation of~\eqref{PDE}, called Reduced-Order Model (ROM). Projection-based methods, such as the Proper Orthogonal Decomposition (POD~\cite{Lumley, Sirovich, Chatterjee}) and the Reduced Basis Method (RBM~\cite{RBmethodPrudhomme}), employ linear data compression (or dimensionality reduction) techniques in order to construct a Reduced-Order Basis (ROB) and use the Galerkin method to derive the ROM. The dimension of the ROB is denoted by $\hat{N}$.  Unlike the HDM, the ROM involves dense reduced operators, as the supports of the reduced-order basis vectors cover the complete domain $\Omega$. Hence, practical speedups, defined as the ratio of the HDM and ROM computation times, are only obtained if $\hat{N}\ll\mathcal{N}$. The reducibility of the parametrized problem is evaluated by inspecting the rate of decay of the Kolmogorov $N$-width $d_{N}(\mathcal{M})$ of the solution manifold $\mathcal{M}$ with respect to the dimension $N$ of an optimal approximation space. We call such spaces optimal $N$-ROM subspaces. More precisely, the Kolmogorov $N$-width quantifies the worst projection error on an optimal $N$-ROM subspace:
\begin{equation*}
d_{N}(\mathcal{M})_{\mathcal{H}} := \underset{\mathcal{H}_{N} \in \textrm{Gr}(N, \mathcal{H})}{\inf} \  \underset{u \in \mathcal{M}}{\sup}  \ \underset{v \in \mathcal{H}_N}{\inf} || u - v ||_{\mathcal{H}} = \underset{\mathcal{H}_{N} \in \textrm{Gr}(N, \mathcal{H})}{\inf} \  \underset{u \in \mathcal{M}}{\sup}  \  || u - \pi_{\mathcal{H}_N}(u) ||_{\mathcal{H}},
\label{KolmogorovWidth}
\end{equation*}
where the Grassmannian $\textrm{Gr}(N, \mathcal{H})$ is the set of all $N$-dimensional subspaces of $\mathcal{H}$ and $\pi_{\mathcal{H}_N}$ is the orthogonal projection onto the subspace $\mathcal{H}_N$. 

\begin{remark}[ROB v.s. optimal $N$-ROM subspace]
The construction of the optimal $N$-ROM subspace is intractable in practice, and a ROB is computed as its approximation by applying an algorithm, such as the POD or the RBM.
Estimating the worst projection error over such a ROB given some assumptions on the Kolmogorov width rate of decay is an active area of research. 
\end{remark}

In this work, $x$ is modeled by a random variable $X$ taking values in $\mathcal{X}$ and following a probability distribution denoted by $p_X : \mathcal{X} \rightarrow \mathbb{R}_+$. 
Another random variable can be defined using the solution application: $U := u(X)$, whose probability distribution $p_U$ is obtained when $X$ follows $p_X$. In particular, the solution manifold is the support of the probability density function $p_U$: $\mathcal{M} = \textrm{supp}(p_U)$. A variant of the Kolmogorov $N$-width is introduced in~\cite{bachmayr2017kolmogorov}, where the worst projection error is replaced by the mean squared error, and is well-suited for this probabilistic context. In this work, we adopt this definition of the Kolmogorov width, but we replace the absolute projection error by the relative projection error, and call it normalized Kolmogorov width:
\begin{equation}
\begin{array}{rcl}
\check{d}_{N}(p_U)_{\mathcal{H}}  & := & \left( \underset{\mathcal{H}_{N} \in \textrm{Gr}(N, \mathcal{H})}{\inf} \  \E_{U \sim p_U} \left[ \eta \left( U, \mathcal{H}_N \right)_{\mathcal{H}}^2 \right] \right)^{1/2} \\
 & = & \left( \underset{\mathcal{H}_{N} \in \textrm{Gr}(N, \mathcal{H})}{\inf} \  \E_{X \sim p_X} \left[ \eta \left( u(X), \mathcal{H}_N \right)_{\mathcal{H}}^2 \right] \right)^{1/2},
\end{array}
\label{NormalizedKWexpected1}
\end{equation}
where $\displaystyle \eta(u,\mathcal{H}_N)_{\mathcal{H}} := \frac{|| u - \pi_{\mathcal{H}_N}(u) ||_{\mathcal{H}}}{|| u ||_{\mathcal{H}}} = \underset{v \in \mathcal{H}_N}{\inf} \frac{|| u - v ||_{\mathcal{H}}}{|| u ||_{\mathcal{H}}}$ is the relative projection error of $u$ onto $\mathcal{H}_N$.

When the solution manifold cannot be embedded in a low-dimensional subspace, the Kolmogorov width decays too slowly with respect to the dimension of the approximation space, and practical speedups cannot be obtained for reasonable accuracies of the reduced predictions: the problem is qualified as non-reducible. Some recent works consider nonlinear data compression and derive optimization problems for the latent dynamics (the evolution of the coefficients of the reduced representation of the solution)~\cite{LEE2020108973, KimChoi2020}. In~\cite{LocalROB, LocalROB2}, the solution manifold $\mathcal{M}$ is partitioned into $K$ subsets $\mathcal{M}_k$ using k-means clustering~\cite{kmeans} in the solution space, which enables constructing a collection of $K$ low-dimensional approximation spaces called dictionary of local ROMs. This paper focuses on dictionaries of local ROMs for non-reducible problems.

\begin{remark}
The relative projection error can be preferred over the absolute error when the norm of the solution changes significantly over the manifold $\mathcal{M}$. Contrary to the absolute error, it does not depend on the magnitude of the solution, and it is symmetric when evaluated between two solutions: $\eta(u, \mathrm{span}(\{v\}))_{\mathcal{H}} = \eta(v,\mathrm{span}(\{u\}))_{\mathcal{H}}$, which enables interpreting it as a dissimilarity measure. Moreover, the reducibility of different problems can be more easily compared via their normalized Kolmogorov widths. It is also common practice to plot the normalized singular values of the POD to manipulate percentages.
\label{Rmk:RelProjErrNormKW}
\end{remark}

\section{Optimal $K$-$N$-ROM-dictionary partitions}
\label{sec:optimalPartitions}

Let $\mathcal{P} \subset \mathcal{M} \subset \mathcal{H} \setminus \{0\}$ be a subset of the solution manifold. The probability $\Prob(\mathcal{P})$ of the event $U \in \mathcal{P}$ reads:
\begin{equation*}
\Prob(\mathcal{P}) = \int_{\mathcal{P}} p_U (u) du = \E_{U \sim p_U} [ \mathds{1}_{\mathcal{P}}(U) ],
\end{equation*} where $\mathds{1}$ is the indicator function. A partition $\{\mathcal{M}_k\}_{1\leq k\leq K}$ of $\mathcal{M}$ is a collection of non-empty subsets of $\mathcal{M}$ such that any point $u$ of the solution manifold $\mathcal{M}$ belongs to exactly one of these subsets. Let $K \geq 2$, and consider a partition of $\mathcal{M}$ into $K$ subsets.
The following definition introduces \textit{optimal $K$-$N$-ROM-dictionary partitions} as partitions that are optimal for reduced-order modeling.
\begin{definition}[Optimal $K$-$N$-ROM-dictionary partitions]
\label{def:optimalPartition}
The family of sets $\{ \mathcal{M}_k \}_{1\leq k \leq K}$ is an optimal $K$-$N$-ROM-dictionary partition of $\mathcal{M}$ if it is a partition of $\mathcal{M}$ and
\begin{equation}
\{ \mathcal{M}_k \}_{1\leq k \leq K} := \mathrm{arg} \underset{ \underset{\mathrm{partition \ of \ }\mathcal{M}}{\{ \mathcal{P}_k \}_{1 \leq k \leq K} }}{\inf} \sum_{k=1}^{K} \Prob(\mathcal{P}_k) \   \check{d}_{N}^2 (p_{U | u \in \mathcal{P}_k})_{\mathcal{H}} .
\label{EqDefOptimalPartitions}
\end{equation}
\end{definition}

An optimal $K$-$N$-ROM-dictionary partition of a solution manifold is a partition of size $K$, leading to a dictionary of local ROMs with $N$ modes minimizing the expectation of the squared intra-cluster normalized Kolmogorov $N$-width. Using Equation~\eqref{NormalizedKWexpected1}, Equation~\eqref{EqDefOptimalPartitions} reads:
\begin{equation}
\{ \mathcal{M}_k \}_{1\leq k \leq K} = \mathrm{arg} \underset{ \underset{\mathrm{partition \ of \ }\mathcal{M}}{\{ \mathcal{P}_k \}_{1 \leq k \leq K} }}{\inf} \sum_{k=1}^{K} \Prob(\mathcal{P}_k) \ \underset{\mathcal{H}_{N}^k \in \textrm{Gr}(N, \mathcal{H})}{\inf} \  \E_{U \sim p_{U | u \in \mathcal{P}_k}} \left[ \eta \left( U, \mathcal{H}_N^k \right)_{\mathcal{H}}^2 \right] .
\label{DetailedCostFctOptPartition}
\end{equation}

\section{Sine dissimilarity}
\label{sec:sineDissimilarity}



We proposed in~\cite{daniel2021physicsinformed} a dissimilarity measure involving the sines of the principal angles between elementary approximation spaces, for time-dependent problems where the solutions can be seen as trajectories in $\mathcal{H}$. 
We referred to it as ROM-oriented dissimilarity and showed that it has the properties of a pseudometric. 
In this work, $\mathcal{M}\subset \mathcal{H}\subset L^2(\Omega)$, so that this sine dissimilarity measure is computed between nonzero elements of $\mathcal{H}$:
\begin{equation}
\tilde{\delta}_{1}(u , v)_{\mathcal{H}} := \sin \measuredangle_{\mathcal{H}} \left(u , v \right) = \sqrt{1 - \frac{ \langle u , v \rangle_{\mathcal{H}}^2}{|| u ||_{\mathcal{H}}^2 || v ||_{\mathcal{H}}^2} },
\label{DefDeltaTildeUn}
\end{equation}
for $(u,v)\in (\mathcal{H} \setminus \{0\} )^2$. Contrary to the distances $\|\cdot\|_{\mathcal{H}}$ and $\|\cdot\|_{L^2(\Omega)}$ commonly used for the construction of dictionaries of local ROMs, the sine dissimilarity measure focuses on the shape of the solutions and is not affected by their intensities. Indeed, for all $(u,v)\in (\mathcal{H} \setminus \{0\} )^2$ and for all $(\lambda_1,\lambda_2) \in \mathbb{R}^{*2}$, $\tilde{\delta}_{1}(\lambda_1 u , \lambda_2 v)_{\mathcal{H}} = \tilde{\delta}_{1}(u , v)_{\mathcal{H}}$. For $(u,v)\in (\mathcal{H} \setminus \{0\} )^2$, let us introduce the following binary relation:
\begin{equation*}
u \sim_{\tilde{\delta}_{1}} v \iff \tilde{\delta}_{1}(u,v)_{\mathcal{H}} = 0.
\end{equation*}
This binary relation is reflexive and symmetric. In addition, it is also transitive, because $\tilde{\delta}_{1}(u , v)_{\mathcal{H}}$ is zero if and only if $u$ and $v$ are linearly dependent, according to the equality case of the Cauchy-Schwarz inequality. This binary relation is thus an equivalence relation, and enables to define the following equivalence classes for the elements of $\displaystyle \mathcal{H}\setminus\{0\}$: $[u] := \{ v \in \mathcal{H}\setminus\{0\} \ | \  v \sim_{\tilde{\delta}_{1}} u  \}$. The quotient set $\mathcal{H}/\sim_{\tilde{\delta}_{1}}$ is defined as the set of all these equivalence classes. In particular in $\mathcal{M}$, collinear solutions are represented by the same element of this quotient set, and $\mathcal{M}/\sim_{\tilde{\delta}_{1}}$ can be seen as the set of directions covered by the solution manifold $\mathcal{M}$. The sine dissimilarity is a metric on $\mathcal{H}/\sim_{\tilde{\delta}_{1}}$.

\begin{proper}
\label{prop:sineToProjError}
The relative projection error can be expressed using the sine dissimilarity: for all $u\in \mathcal{H} \setminus \{0\}$ and $\mathcal{H}_{N} \in \mathrm{Gr}(N, \mathcal{H})$,
\begin{equation*}
\eta \left( u, \mathcal{H}_{N} \right)_{\mathcal{H}}^2 = 1 - N + \sum_{j=1}^N \tilde{\delta}_{1}(u, h_{j})_{\mathcal{H}}^2,
\end{equation*} 
where the $N$ functions $h_{j}$ form an orthonormal basis of the subspace $\mathcal{H}_{N}$.
\end{proper}
\begin{proof}
The orthogonal projection of $u\in\mathcal{H}\setminus\{0\}$ onto the subspace $\mathcal{H}_N$ reads:
$\pi_{\mathcal{H}_N}(u) = \sum_{j=1}^N \langle u , h_j \rangle_{\mathcal{H}} h_j $.
Therefore:
\begin{equation*}
\begin{array}{rcl}
\eta \left( u, \mathcal{H}_{N} \right)_{\mathcal{H}}^2 
 & = & || u ||_{\mathcal{H}}^{-2} \left( || u ||_{\mathcal{H}}^2 - 2 \displaystyle\sum_{j=1}^N \langle u , h_j \rangle_{\mathcal{H}}^2 + \displaystyle\sum_{j=1}^N \sum_{i=1}^N \langle u , h_j \rangle_{\mathcal{H}} \langle u , h_i \rangle_{\mathcal{H}} \langle h_i , h_j \rangle_{\mathcal{H}} \right) \\
 & = & 1 - \displaystyle\sum_{j=1}^N \frac{\langle u , h_j \rangle_{\mathcal{H}}^2}{|| u ||_{\mathcal{H}}^{2}}.
\end{array}
\end{equation*}
The proof is ended using the orthonormality of the basis $\{ h_j \}_{1\leq j \leq N}$ of $\mathcal{H}_N$ and Equation~\eqref{DefDeltaTildeUn}.
\end{proof}

\begin{cor}
\label{1DProjErrorVSdissimilarity}
For all $(u,v)\in (\mathcal{H} \setminus \{0\} )^2$,
\begin{equation*}
\tilde{\delta}_{1}(u , v)_{\mathcal{H}} = \eta \left( u, \mathrm{span}(\{v\}) \right)_{\mathcal{H}}.
\end{equation*}
\end{cor}


\section{Optimal $K$-$1$-ROM-dictionary partitions}

Taking $N=1$, Equation~\eqref{DetailedCostFctOptPartition} defines optimal $K$-$1$-ROM-dictionary partitions as the solutions of the optimization problem:
\begin{equation}
\{ \mathcal{M}_k \}_{1\leq k \leq K} = \mathrm{arg} \underset{ \underset{\mathrm{partition \ of \ }\mathcal{M}}{\{ \mathcal{P}_k \}_{1 \leq k \leq K} }}{\inf} \sum_{k=1}^{K} \Prob(\mathcal{P}_k) \ \underset{\tilde{u}_k \in \mathcal{H}/\sim_{\tilde{\delta}_1}}{\inf} \  \E_{U \sim p_{U | u \in \mathcal{P}_k}} \left[ \tilde{\delta}_{1}(U , \tilde{u}_k)_{\mathcal{H}}^2 \right] ,
\label{CostFctN1Continuous}
\end{equation}
where the relative projection error is replaced by the sine dissimilarity thanks to Corollary~\ref{1DProjErrorVSdissimilarity}. Note that the second infimum is taken on the quotient set $\mathcal{H}/\sim_{\tilde{\delta}_1}$, because looking for an optimal element of the quotient set is equivalent to searching for an optimal 1D approximation space in $\mathrm{Gr}(1,\mathcal{H})$. Equation~\eqref{CostFctN1Continuous} can be interpreted as the continuous version of a representative-based clustering problem, where for a given integer $K \geq 2$, the objective is to find $K$ representative elements whose nearest neighbors for a given dissimilarity measure define the $K$ clusters. Given a metric measure space $\mathcal{V}$ whose metric (\textit{resp.} measure) is denoted by $\delta_{\mathcal{V}}$ (\textit{resp.} $\mu_{\mathcal{V}}$), and given a subset $\mathcal{V'}$ of $\mathcal{V}$ with a nonzero measure, the continuous representative-based clustering problem can be stated as follows:
\begin{equation*}
\underset{ \underset{\mathrm{partition \ of \ }\mathcal{V}'}{\{ \mathcal{P}_k \}_{1 \leq k \leq K} }}{\inf} \sum_{k=1}^{K} \underset{\tilde{u}_k \in \mathcal{V}}{\inf} \  \int_{\mathcal{P}_k} \delta_{\mathcal{V}}^2 (v,\tilde{u}_k) d\mu_{\mathcal{V}}(v),
\end{equation*} which is equivalent to:
\begin{equation*}
\underset{ \underset{\mathrm{partition \ of \ }\mathcal{V}'}{\{ \mathcal{P}_k \}_{1 \leq k \leq K} }}{\inf} \sum_{k=1}^{K} \frac{\mu_{\mathcal{V}}(\mathcal{P}_k)}{\mu_{\mathcal{V}} (\mathcal{V}')} \underset{\tilde{u}_k \in \mathcal{V}}{\inf} \  \frac{1}{\mu_{\mathcal{V}}(\mathcal{P}_k)} \int_{\mathcal{P}_k} \delta_{\mathcal{V}}^2 (v,\tilde{u}_k) d \mu_{\mathcal{V}}(v).
\end{equation*}
The ratio $\mu_{\mathcal{V}}(\mathcal{P}_k)/\mu_{\mathcal{V}} (\mathcal{V}')$ can be seen as the probability $\Prob(\mathcal{P}_k)$ of being in $\mathcal{P}_k$ when drawing a realization of the uniform distribution on $\mathcal{V}'$. The integral term normalized by the measure of the cluster corresponds to the expectation in Equation~\eqref{CostFctN1Continuous}. In both cases, the integrand is a squared metric. When $\mathcal{V}$ is a Hilbert space and $\delta_{\mathcal{V}}$ is the norm induced by its inner product, the aforementioned clustering problem is a continuous k-means clustering problem, where the objective is to find clusters minimizing the sum of the intra-cluster inertia. In this case, the optimal representative elements are the centroids (or means) of the clusters. In a nutshell, the optimal $K$-$1$-ROM-dictionary partitions are the solutions of a representative-based clustering problem on the quotient set $\mathcal{H}/\sim_{\tilde{\delta}_1}$ endowed with the metric $\tilde{\delta}_1$. 

\section{Algorithm for the construction of a dictionary of local ROMs}
\label{sec:algo}

The optimization problem defining optimal $K$-$N$-ROM-dictionary partitions in Definition~\ref{def:optimalPartition} is numerically intractable, since there is an infinite number of possible partitions of the continuous solution manifold and since each candidate partition requires solving $K$ optimization problems for the construction of the local ROMs.
In the same fashion the (snapshot)-POD has been proposed as a practical procedure for approximating the optimal $N$-ROM subspace, we propose an algorithm approximating the optimal $K$-$N$-ROM-dictionary partitions, given a sampling of the solution manifold as a set of $m$ precomputed solutions $\widehat{\mathcal{M}} := \{ u_i \}_{1 \leq i \leq m}$. This \emph{a priori} sample is usually done by applying a design of experiments over the parameter domain~$\mathcal{X}$. 

\subsection{Approximate optimal $K$-$N$-ROM-dictionary partitions of discrete solution sets}

We look for an approximation of the optimal $K$-$N$-ROM-dictionary partitions $\{ \widehat{\mathcal{M}}_k \}_{1 \leq k \leq K}$ of the discretized set $\widehat{\mathcal{M}}$.
The probabilities $\Prob(\widehat{\mathcal{M}}_k)$ are obtained by taking the ratios $| \widehat{\mathcal{M}}_k | / m$, where $| \widehat{\mathcal{M}}_k |$ is the cardinality of $\widehat{\mathcal{M}}_k$.
The true but unknown probability density functions $p_{U | u \in \mathcal{M}_k}$ are replaced by the empirical distributions
$\displaystyle \widehat{p}_{U | u \in \widehat{\mathcal{M}}_k} (u) = \frac{1}{|\widehat{\mathcal{M}}_k |} \sum_{i=1}^{m} \mathds{1}_{\widehat{\mathcal{M}}_k}(u_i) \ \delta(u - u_i)$,
where $\delta$ is the Dirac delta function.
Using Equation~\eqref{NormalizedKWexpected1}, the squared normalized Kolmogorov $N$-width of this probability mass function is then:
\begin{equation*}
\check{d}_{N}(\widehat{p}_{U | u \in \widehat{\mathcal{M}}_k})_\mathcal{H}^2 = \underset{\mathcal{H}_{N}^k \in \textrm{Gr}(N, \mathcal{H})}{\inf} \  \frac{1}{|\widehat{\mathcal{M}}_k |} \sum_{i=1}^{m}\mathds{1}_{\widehat{\mathcal{M}}_k}(u_i) \ \eta \left( u_i, \mathcal{H}_N \right)_{\mathcal{H}}^2.
\end{equation*}
Like the RBM and the snapshot-POD methods, and to derive a computable algorithm, a basis for the approximation of the best subspace $\mathcal{H}_{N}^k$ is searched in the set  $\mathcal{A}_{N}(\widehat{\mathcal{M}}_k)$, $N \leq \dim(\textrm{span}(\widehat{\mathcal{M}}_k))$, defined as the set containing all the $\mathcal{H}$-orthonormal families of $N$ elements of $\mathrm{span}(\widehat{\mathcal{M}}_k)$.
From Definition~\ref{def:optimalPartition} and Property~\ref{prop:sineToProjError}, an approximation of the optimal $K$-$N$-ROM-dictionary partitions $\{ \widehat{\mathcal{M}}_k \}_{1 \leq k \leq K}$ of a discrete solution set $\widehat{\mathcal{M}}$ is sought as:
\begin{equation}
\label{eq:approxPartitions}
{\rm arg} \ \underset{\underset{\mathrm{partition \ of \ } \widehat{\mathcal{M}}}{ \{ \widehat{\mathcal{P}}_k \}_{1\leq k \leq K}}}{\min} \ \sum_{k=1}^{K} \ \underset{{ \{ h_j^k\}_{1\leq j\leq N}\in \mathcal{A}_{N}(\widehat{\mathcal{M}})}}{\min} \  \sum_{i=1}^{m} \mathds{1}_{\widehat{\mathcal{P}}_k}(u_i) \  \sum_{j=1}^N \tilde{\delta}_{1}(u_i, h^k_{j})_{\mathcal{H}}^2.
\end{equation}

\subsection{Approximate optimal $K$-$1$-ROM-dictionary partitions of discrete solution sets}

\begin{proper}
\label{prop:equivalence}
When considering a discrete solution set $\widehat{\mathcal{M}}$ and if the 1-ROM subspaces are sought in $\mathrm{span}(\widehat{\mathcal{M}})$, the optimal $K$-$1$-ROM-dictionary partitions are exactly the minimizers of the cost function of k-medoids clustering with the sine dissimilarity measure $\tilde{\delta}_{1}$.
\end{proper}
\begin{proof}
From Equation~\eqref{eq:approxPartitions} (or discretizing the solution manifold in Equation~\eqref{CostFctN1Continuous}), optimal $K$-$1$-ROM-dictionary partitions satisfy
\begin{equation*}
{\rm arg} \ \underset{\underset{\mathrm{partition \ of \ } \widehat{\mathcal{M}}}{ \{ \widehat{\mathcal{P}}_k \}_{1\leq k \leq K}}}{\min} \ \sum_{k=1}^{K} \ \underset{\tilde{u}_k \in \widehat{\mathcal{M}}}{\min} \  \sum_{i=1}^{m} \mathds{1}_{\widehat{\mathcal{P}}_k}(u_i) \  \tilde{\delta}_{1}(u_i, \tilde{u}_k)_{\mathcal{H}}^2,
\end{equation*}
from which we recognize the cost function of k-medoids clustering.
\end{proof}

Hence, the case $N=1$ leads to an optimization problem for which various computable heuristic approaches have been proposed, including the Partitioning Around Medoids (PAM~\cite{kMedoidsPAM}).

\subsection{Algorithm for approximate optimal partitions}

Property~\ref{prop:equivalence} cannot be directly extended to $N>1$, and optimizing the partition and the approximation spaces simultaneously in Equation~\eqref{eq:approxPartitions} requires computing many candidate local subspaces, which is very expensive. As a practical algorithm, we propose to: (i) compute the optimal $K$-$1$-ROM-dictionary partitions over the sampled solution manifold, and (ii) compute the local $N$-dimensional approximation spaces using any classical reduced-order modeling method on each element of this partition, for instance the snapshot-POD or the RBM.

We applied this algorithm in our previous work~\cite{daniel2021physicsinformed, UQindustrialDesign}, justifying the use of the sine dissimilarity by qualitative arguments saying that an angle-based dissimilarity should be preferred when searching for approximating vector subspaces. In the present work, we introduce a natural concept of optimal partition of the solution manifold for reduced-order modeling, and obtain this algorithm as a computable approximation of this partition. In the literature, dictionaries of ROMs have been proposed by partitioning the solution manifold with the Euclidean or the $\|\cdot\|_{L^2(\Omega)}$ distances~\cite{LocalROB, LocalROB2}, or by partitioning the parameter domain $\mathcal{X}$~\cite{parmSpacePartitioning}. For dictionaries of local linear ROMs, we show that a representative-based clustering procedure involving the sine dissimilarity should be used in order to get optimal partitions in terms of normalized Kolmogorov widths.

\begin{remark}
In addition to the arguments given in Remark~\ref{Rmk:RelProjErrNormKW} for the use of normalized Kolmogorov widths, using relative errors enables linking the concept of optimal partitions with a representative-based cluster analysis. Trying to do the same for absolute Kolmogorov widths would have required considering a dissimilarity obtained by symmetrizing the absolute errors: the link between the corresponding optimal $K$-$1$-ROM-dictionary partitions and such dissimilary would be lost, and Property~\ref{prop:equivalence} would not hold anymore.
\end{remark}

\section{Application}


This procedure for the construction of local ROMs has been applied to a real industrial problem taken from the aeronautical industry, where the pPDE describes the mechanical behavior of a high-pressure turbine blade in an aircraft engine, subjected to uncertain thermal loadings modeled by random fields. For more details about this industrial application, see~\cite{UQindustrialDesign}.

\section{Outlook}

Section~\ref{sec:algo} started by considering an \emph{a priori} sampling of the solution manifold.
We can introduce the problem of simultaneously optimizing the partition of the solution manifold, the construction of the local approximation spaces, and the solution manifold sampling. In the same fashion as the Reduced Basis Method, this means being able to compute efficient and sharp \emph{a posteriori} relative error bounds.

\section*{Acknowledgements}
\noindent Study funded by Safran and ANRT (Association Nationale de la Recherche et de la Technologie, France).

\appendix

\bibliographystyle{unsrt}
\bibliography{Biblio}

\end{document}